\documentclass[preprint,12pt]{elsarticle}

\usepackage{graphicx}
\usepackage{amssymb}
\usepackage{url}
\usepackage{amsmath}
\usepackage{siunitx}
\usepackage{graphicx} 
\usepackage{float} 
\usepackage{subfigure} 

\usepackage{lineno}

\journal{arxiv}

\begin{document}

\begin{frontmatter}


\title{Existence and Uniqueness of Reynolds equation under natural boundary conditions}

\author{Wang Qun}

\begin{abstract}
We often use the finite element calculation to solve the Reynolds equation(a kind of elliptic differential equations) under the natural boundary conditions.During numerical solution,we found an interesting phenomenon : if we start the calculate from all pressure is 0,the area of pressure greater than 0 is always 'growing' as the iteration progresses. Inspired by this phenomenon,we prove the existence and uniqueness of Reynolds equation under natural boundary conditions(double Reynolds boundary condition).
\\\\
We can solve the Reynolds equation under given boundary conditions which set the pressure is greater
than 0 in the boundary $\Omega$ and the pressure is zero on the boundary $\Omega$ :
\[\begin{array}{l}
p(x,y) = 0,(x,y) \in \Omega \\
p(x,y) > 0,(x,y) \in {\Omega _{in}}
\end{array}\]
Record $S$ as the set of this kind of boundary conditions.The partial order relations and addition operation established on $S$.Record $A,B$ is the boundary conditions and $\Omega ^A,\Omega ^B$ is the corresponding boundary:
\[\begin{array}{c}
A \ge B \leftrightarrow \left[ {(x,y) \in \Omega _{in}^B \to (x,y) \in \Omega _{in}^A} \right]\\
C=A+B  \leftrightarrow \left[ {(x,y) \in \Omega _{in}^A \cup (x,y) \in \Omega _{in}^B \leftrightarrow (x,y) \in \Omega _{in}^{C}} \right]
\end{array}\]
Via the Hopf maximum principle and the linearity of elliptic differential equations,we can get the relationship between the pressure solution under boundary condition A+B and the solution under boundary condition $A$ and $B$, we can get:
\[{p_{A + B}} \ge \max ({p_A},{p_B}) \ge 0\]
so we can get the set $S$ is closed under addition operation:
\[A,B \in S \to A + B \in S\]
By similar method,we can establish an equivalence relation between the natural boundary conditions and the maximum element of S.
\[\neg \exists s \in S,s > A;A \in S \leftrightarrow A\;meets\;natural\;boundary\;conditions\]
As the boundary of the whole calculate domain ${\overline D }$ must be a upper bound of the set S,so there must be a greatest element $G\left( S \right)$ in set S:
\[\left\{ {\begin{array}{*{20}{l}}
{\forall s \in S,{\overline D } \ge s}\\
{A,B \in S \to A + B \in S}
\end{array}} \right. \to \exists G\left( S \right) \in S\left( {\forall C \in S,C \le G\left( S \right)} \right)\]
Using the above equivalence relation,the existence and uniqueness of Reynolds equation under natural boundary conditions can be proved.\\\\
We give a method to solve the one dimensional numerical Reynolds equation,which has the complexity of O(N).As the discretization of Reynolds equation also has the complexity of O(N),the overall solution complexity will not be less than O(N). So this method's speed reach the theoretical upper bound meanwhile has no convergence error.We publish our code on Github.
\\\\
We discuss the physical meaning of our paper in the last of paper.
\\
\end{abstract}

\begin{keyword}
Reynolds equation \sep Elliptic Differential Equations \sep Existence and Uniqueness \sep numerical method \sep Hopf maximum principle \sep Reynolds boundary condition

\end{keyword}

\end{frontmatter}


\section{Introduction}
The study of free boundary problems in the hydrodynamic lubrication gave place to many
works covering some fundamental and applied aspects\cite{bayada2007existence}\cite{rapidmethod}. For lubricated devices with Newtonian fluid, the classical thin film Reynolds equation is\cite{wen2012principles}
\[\frac{\partial }{{\partial x}}\left( {\frac{{\rho {h^2}}}{\eta }\frac{{\partial p}}{{\partial x}}} \right) + \frac{\partial }{{\partial y}}\left( {\frac{{\rho {h^2}}}{\eta }\frac{{\partial p}}{{\partial y}}} \right) = 6\left( {{V_x}\frac{{\partial \rho h}}{{\partial x}} + {V_y}\frac{{\partial \rho h}}{{\partial y}} + 2\frac{{\partial \rho h}}{{\partial t}}} \right)\]
where $p$ is the pressure of the lubricant, $h$ is the given film thickness, $\rho$ is the density of
lubricant, $\eta$ is the viscosity of the lubricant, $V_x,V_y$ is the relative velocity in the $x,y$ direction of the surfaces in which lubricant takes place and $t$ is the time.

As $h$ , $\rho$ , $\eta$ , $V_x,V_y$ and $ \frac {{\partial \rho h}}{{\partial t}}$ is given by $x,y,t$ ,so we can simplify the equation to 
\begin{equation} \label{ori}
\frac{\partial }{{\partial x}}\left( {{f_1}\frac{{\partial p}}{{\partial x}}} \right) + \frac{\partial }{{\partial y}}\left( {{f_2}\frac{{\partial p}}{{\partial y}}} \right) = {f_3}
\end{equation} 
where $f_1,f_2\ge 0$.

Equation \ref{ori} is to be used only on an unknown part of $\Omega_{in}$ in which pressure is strictly greater
than the vaporization pressure (often taken as zero). In the other part, namely cavitation region,
Reynolds equation is no longer valid and the pressure is the vaporization pressure.So the natural boundary conditions(or named as Reynolds boundary condition) can be note as :
\begin{equation}\label{condition}
\left\{ \begin{array}{l}
{p(x,y)} > 0,(x,y) \in {\Omega _{in}}\\
{f_3}(x,y) >  = 0,(x,y) \notin {\Omega _{in}}\\
p\left| {_\Omega } \right. = 0\\
\frac{{\partial p}}{{\partial \Omega_n }} =  0
\end{array} \right.
\end{equation}

One of The key point of this paper is defining a series of operations on the set of boundary conditions,which is in part \ref{condition}.The other one is the set of boundary condition with positive solution is closed under addition operation,which is in part \ref{close}.We give an example in part \ref{example} to explain part \ref{condition} and part \ref{close} .The equivalence of our define and natural boundary condition is in part \ref{equivalence}.

For engineers, we often use a zero setting method to find the solution of Reynolds equation. We can use the method in this paper to prove that the method has and only has one solution with iteration error of 0 in part \ref{numerical_proof}. A fast numerical method for one dimensional is given in part \ref{fast}.

We give conclusions in part \ref{Conclusions}.

For mathematician,part \ref{condition},\ref{close},\ref{Conclusions} is important;part 
 \ref{example},\ref{equivalence} is auxiliary and trivial.
 
For engineer,part \ref{numerical_proof},\ref{fast} is important,part \ref{condition} is for reference.

\section{Addition operation and partial order relation of boundary conditions}\label{opera}
The key point to solve Reynolds equation under natural boundary conditions is to find the unknown boundary of the fully lubricated area $\Omega$ in which pressure is strictly greater than zero.
A method to find $\Omega$ is to assume that $\Omega$ equals to 
$\Omega_{1}$,and solve the Reynolds equation under the assuming boundary conditions:
\[{p_1}(x,y) = 0,(x,y) \in {\Omega _1}\]
If $p_1$ and $\Omega_{1}$ can meet formula \ref{condition},the solution is finded.If not, we can adjust the boundary according to the pressure obtained to find the boundary that meets the conditions

In this way, all possible boundary conditions and  can form a set,and each element has a corresponding solution and boundary of pressure greater than 0. All the boundary conditions that make the pressure of the solution positive can form a subset $S$. This paper focuses on the partial order relation in set S and the relation between addition operation and natural boundary.

For more accurate description and proof,we define the partial order relation and addition operation on the assuming boundary conditions:
\begin{equation}\label{operation}
\begin{array}{c}
A \ge B \leftrightarrow \left[ {(x,y) \in \Omega _{in}^B \to (x,y) \in \Omega _{in}^A} \right]\\
C=A+B  \leftrightarrow \left[ {(x,y) \in \Omega _{in}^A \cup (x,y) \in \Omega _{in}^B \leftrightarrow (x,y) \in \Omega _{in}^{C}} \right]
\end{array}
\end{equation}
which can be easily understand as:
\begin{itemize}
\item Boundary condition A is not less than B means all points in $\Omega _{in}^B$ are also in $\Omega _{in}^A$
\item Boundary condition C is the sum of A and B means all points in $\Omega _{in}^A$ or $\Omega _{in}^B$ are also in $\Omega _{in}^C$ and all points in $\Omega _{in}^C$ is also in $\Omega _{in}^A$ or $\Omega _{in}^B$.
\end{itemize}
What needs to be reminded is the difference between partial order relation and total order relation:A and B can have no order relationship like A is $\{p(x)>0,-1<x<1\}$ and B is $\{p(x)>0,0<x<2\}$.So there is difference in maximum and greatest: A is maximum means no one element is greater than A;A is greatest means A is greater than every one element else.

\section{The set of boundary condition is closed under addition operation}\label{close}
Take any two boundary condition A and B from the set S, record the solution under boundary condition A and B as $p_A$ and $p_B$.By the properties of set s:
\[\left\{ {\begin{array}{*{20}{l}}
{{p_A}(x,y) > 0,(x,y) \in \Omega _{_{in}}^A;{p_A}(x,y) = 0,(x,y) \in \Omega _{_{in}}^{A + B} \cap (x,y) \notin \Omega _{_{in}}^A}\\
{{p_B}(x,y) > 0,(x,y) \in \Omega _{_{in}}^B;{p_B}(x,y) = 0,(x,y) \in \Omega _{_{in}}^{A + B} \cap (x,y) \notin \Omega _{_{in}}^B}
\end{array}} \right.\]
Record $p_m$ as the max of $p_A,p_B$,and the solution under boundary condition of $A+B$ as $p_{A+B}$.Record $p_m^\Delta$ as the difference of $p_{A+B}$ and $p_m$:
\begin{equation}\label{delta}
p_{A+B}=p_m+p_m^\Delta
\end{equation}
We will prove that:
\begin{equation}\label{larger}
\frac{\partial }{{\partial x}}\left( {{f_1}\frac{{\partial {p_m}}}{{\partial x}}} \right) + \frac{\partial }{{\partial y}}\left( {{f_2}\frac{{\partial {p_m}}}{{\partial y}}} \right) - {f_3} \ge 0,\left( {x,y} \right) \in \Omega _{in}^{A + B}\end{equation}

According to the above physical definition, the Reynolds equation is strictly satisfied when the pressure is greater than 0.And according to the definition of the addition operation of boundary conditions,at least one of ${p_A}(x,y) > 0,{p_B}(x,y) > 0$ holds.Thirdly,$p_m$ is the max of $p_A,p_B$.We can get that at least one of the following two equations holds:
\begin{equation}\label{both}
{\left\{ {\begin{array}{*{20}{l}}
{{p_m}(x,y) = {p_A}(x,y) > 0 \to \frac{\partial }{{\partial x}}\left( {{f_1}\frac{{\partial {p_A}}}{{\partial x}}} \right) + \frac{\partial }{{\partial y}}\left( {{f_2}\frac{{\partial {p_A}}}{{\partial y}}} \right) - {f_3} = 0}\\
{{p_m}(x,y) = {p_B}(x,y) > 0 \to \frac{\partial }{{\partial x}}\left( {{f_1}\frac{{\partial {p_B}}}{{\partial x}}} \right) + \frac{\partial }{{\partial y}}\left( {{f_2}\frac{{\partial {p_B}}}{{\partial y}}} \right) - {f_3} = 0}
\end{array}} \right.}
\end{equation}
Write the formula \ref{larger},\ref{both} as limit form:
\[\begin{array}{l}
\frac{\partial }{{\partial x}}\left( {{f_1}\frac{{\partial p}}{{\partial x}}} \right) + \frac{\partial }{{\partial y}}\left( {{f_2}\frac{{\partial p}}{{\partial y}}} \right) - {f_3}\\
 = \mathop {\lim }\limits_{\varepsilon  \to 0} {{\left[ {p(x + {f_1}\varepsilon ,y) + p(x - {f_1}\varepsilon ,y) + p(x,y + {f_2}\varepsilon ) + p(x,y - {f_2}\varepsilon ) - 4p(x,y)} \right]} \mathord{\left/
 {\vphantom {{\left[ {p(x + {f_1}\varepsilon ,y) + p(x - {f_1}\varepsilon ,y) + p(x,y + {f_2}\varepsilon ) + p(x,y - {f_2}\varepsilon ) - 4p(x,y)} \right]} \varepsilon }} \right.
 \kern-\nulldelimiterspace} \varepsilon } - {f_3}
\end{array}\]
note that:
\[\begin{array}{l}
{p^{(0)}} = p(x,y)\\
{p^{(1)}} = p(x + {f_1}\varepsilon ,y)\\
{p^{(2)}} = p(x - {f_1}\varepsilon ,y)\\
{p^{(3)}} = p(x,y + {f_2}\varepsilon )\\
{p^{(4)}} = p(x,y - {f_2}\varepsilon )
\end{array}\]
therefore:
\[\begin{array}{l}
{p_m}^{(0)} = {p_A}^{(0)}\\
 \to {p_A}^{(0)} > 0\\
 \to \mathop {\lim }\limits_{\varepsilon  \to 0} \frac{{{p_A}^{(1)} + {p_A}^{(2)} + {p_A}^{(3)} + {p_A}^{(4)} - 4{p_A}^{(0)}}}{\varepsilon } - {f_3} = 0\\
 \to \mathop {\lim }\limits_{\varepsilon  \to 0} \frac{{{p_A}^{(1)} + {p_A}^{(2)} + {p_A}^{(3)} + {p_A}^{(4)} - 4{p_m}^{(0)}}}{\varepsilon } - {f_3} = 0\\
 \to \mathop {\lim }\limits_{\varepsilon  \to 0} \frac{{{p_m}^{(1)} + {p_m}^{(2)} + {p_m}^{(3)} + {p_m}^{(4)} - 4{p_m}^{(0)}}}{\varepsilon } - {f_3} \ge \mathop {\lim }\limits_{\varepsilon  \to 0} \frac{{{p_A}^{(1)} + {p_A}^{(2)} + {p_A}^{(3)} + {p_A}^{(4)} - 4{p_m}^{(0)}}}{\varepsilon } - {f_3} \ge 0
\end{array}\]
similarly:
\[{p_m}^{(0)} = {p_B}^{(0)} \to \mathop {\lim }\limits_{\varepsilon  \to 0} \frac{{{p_m}^{(1)} + {p_m}^{(2)} + {p_m}^{(3)} + {p_m}^{(4)} - 4{p_m}^{(0)}}}{\varepsilon } - {f_3} \ge 0\]
as ${p_m} = \max \left( {{p_A},{p_B}} \right)$,we can get 
\[{p_m} = \max \left( {{p_A},{p_B}} \right) \to \left[ {{p_m}^{(0)} = {p_A}^{(0)}} \right]OR\left[ {{p_m}^{(0)} = {p_B}^{(0)}} \right]\]
Either of above formula must be true, and the result of formula \ref{larger} is obtained,so formula \ref{larger} is proved.

We can also regard $p_m^\Delta$ as a solutions under the boundary condition of $A+B$.Via the linearity of elliptic differential equations,we can get that:
\[\left\{ \begin{array}{l}
p = p_m^\Delta  + {p_m}\\
\frac{\partial }{{\partial x}}\left( {{f_1}\frac{{\partial p}}{{\partial x}}} \right) + \frac{\partial }{{\partial y}}\left( {{f_2}\frac{{\partial p}}{{\partial y}}} \right) - {f_3} = 0\\
\frac{\partial }{{\partial x}}\left( {{f_1}\frac{{\partial {p_m}}}{{\partial x}}} \right) + \frac{\partial }{{\partial y}}\left( {{f_2}\frac{{\partial {p_m}}}{{\partial y}}} \right) - {f_3} \ge 0
\end{array} \right. \to \frac{\partial }{{\partial x}}\left( {{f_1}\frac{{\partial p_m^\Delta }}{{\partial x}}} \right) + \frac{\partial }{{\partial y}}\left( {{f_2}\frac{{\partial p_m^\Delta }}{{\partial y}}} \right) \le 0\]
via the Hopf maximum principle:
\[\frac{\partial }{{\partial x}}\left( {{f_1}\frac{{\partial p_m^\Delta }}{{\partial x}}} \right) + \frac{\partial }{{\partial y}}\left( {{f_2}\frac{{\partial p_m^\Delta }}{{\partial y}}} \right) \le 0 \to p_m^\Delta  \ge 0 \to {p_{A + B}} = {p_m} + p_m^\Delta  \ge 0\]
which means the solution under boundary condition A+B also is positive.So we have proved the set S is closed under under addition operation.
\[A,B \in S \to A + B \in S\]

\section{An example for part \ref{condition} and part \ref{close}}
\label{example}
Take an one dimensional example,considering of the solution of:
\[\left\{ \begin{array}{l}
\frac{{{\partial ^2}p}}{{\partial {x^2}}} =  - 2,( - 3 \ge x \ge 3)\\
p = 0(x =  - 3,x = 3)
\end{array} \right.\]
The solution is :
\[p(x)=9-{x^2}( - 3 \ge x \ge 3)\]
Assume that boundary condition $A$ and $B$ is 
\[\begin{array}{l}
A:p(x) = 0,(x =  - 3,x = 1);p(x) \ge 0,( - 3 < x < 1)\\
B:p(x) = 0,(x =  - 1,x = 3);p(x) \ge 0,( - 1 < x < 3)
\end{array}\]
so the solution of under $A$ and $B$ is :
\[\begin{array}{l}
{p_A}(x) = 4 - {(x + 1)^2}\\
{p_B}(x) = 4 - {(x - 1)^2}
\end{array}\]
as ${p_m} = \max \left( {{p_A},{p_B}} \right)$,we can get:
\[{p_m} = \left\{ \begin{array}{l}
4 - {(x + 1)^2}( - 3 < x \le 0)\\
4 - {(x - 1)^2}(0 < x < 3)
\end{array} \right.\]
Since there are non-differentiable points, we need to introduce Dirac delta function.Record $\delta$ as Dirac delta function,we can get:
\[\left\{ {\begin{array}{*{20}{l}}
{\frac{{{\partial ^2}{p_m}}}{{\partial {x^2}}} =  - 2,( - 3 < x < 0 \cup 0 < x < 3)}\\
{{{\left. {\frac{{{\partial ^2}{p_m}}}{{\partial {x^2}}}} \right|}_{x = 0}} = \delta ({{\left. {\frac{{\partial {p_m}}}{{\partial x}}} \right|}_{x = {0^ + }}} - {{\left. {\frac{{\partial {p_m}}}{{\partial x}}} \right|}_{x = {0^ - }}}) = 4\delta }
\end{array}} \right.\]
According to formula \ref{delta},we can find that
\[\frac{{{\partial ^2}{p_m}^\Delta }}{{\partial {x^2}}} = \frac{{{\partial ^2}p}}{{\partial {x^2}}} - \frac{{{\partial ^2}{p_m}}}{{\partial {x^2}}} \to \left\{ {\begin{array}{*{20}{l}}
{\frac{{{\partial ^2}{p_m}^\Delta }}{{\partial {x^2}}} = 0,( - 3 < x < 0 \cup 0 < x < 3)}\\
{{{\left. {\frac{{{\partial ^2}{p_m}^\Delta }}{{\partial {x^2}}}} \right|}_{x = 0}} =  - 4\delta }
\end{array}} \right.\]
solve it we can get that:
\[{p_m}^\Delta  = \left\{ \begin{array}{l}
2x + 6( - 3 < x \le 0)\\
 - 2x + 6(0 < x < 3)
\end{array} \right.\]
we can verify that
\[{p_m}^\Delta  + {p_m} = \left\{ \begin{array}{l}
2x + 6 + 4 - {(x + 1)^2} = 9 - {x^2}( - 3 < x \le 0)\\
 - 2x + 6 + 4 - {(x - 1)^2} = 9 - {x^2}(0 < x < 3)
\end{array} \right. = p\]
\begin{figure}[H] 
\centering 
\includegraphics[height=0.7\textwidth]{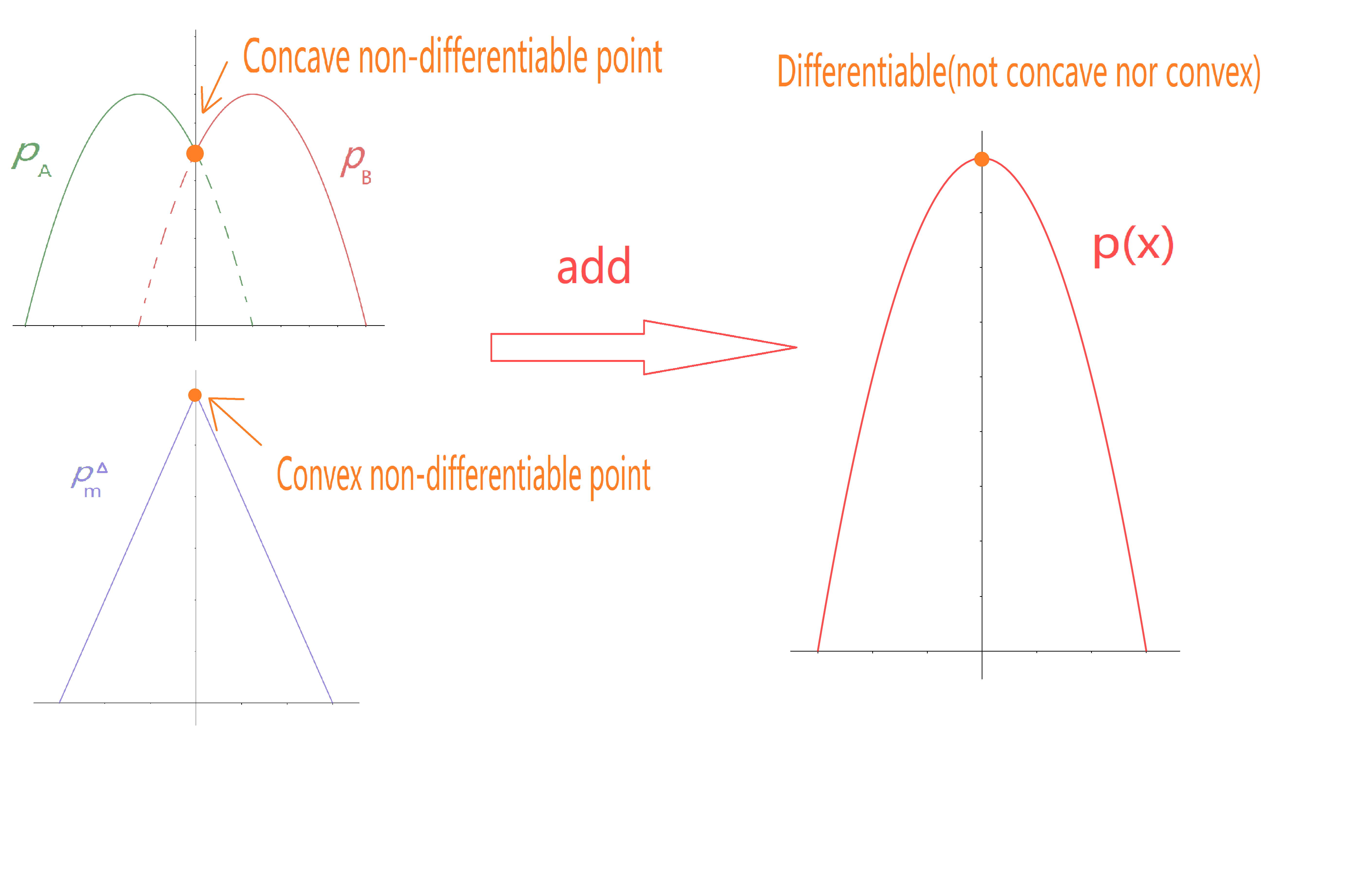} 
\caption{add of domain} 
\label{Fig.main2} 
\end{figure}

\section{The equivalence relation between the natural boundary conditions and the maximum element of S}
\label{equivalence}

\noindent\textbf{proof 5.1}

First we prove boundary condition $A$ meets the natural boundary condition (formula \ref{condition}) must be the maximum element of S. 

Suppose there is a element $B$ is larger than $A$ in $S$:
\[\exists B \in S,B > A\]
define that:
\[{p_B}\left( {x,y} \right) - {p_A}\left( {x,y} \right) = {p_{B - A}}\left( {x,y} \right)\]
 according to 
\[\left\{ \begin{array}{l}
\frac{\partial }{{\partial x}}\left( {{f_1}\frac{{\partial {p_A}}}{{\partial x}}} \right) + \frac{\partial }{{\partial y}}\left( {{f_2}\frac{{\partial {p_A}}}{{\partial y}}} \right) - {f_3} = 0,(x,y) \in \Omega _{in}^A\\
\frac{\partial }{{\partial x}}\left( {{f_1}\frac{{\partial {p_B}}}{{\partial x}}} \right) + \frac{\partial }{{\partial y}}\left( {{f_2}\frac{{\partial {p_B}}}{{\partial y}}} \right) - {f_3} = 0,(x,y) \in \Omega _{in}^B\\
{p_A} = 0,(x,y) \notin \Omega _{in}^A
\end{array} \right.\]
we can find that:
\[{\left\{ \begin{array}{l}
{p_{B - A}} = 0,(x,y) \in {\Omega ^B}\\
\frac{\partial }{{\partial x}}\left( {{f_1}\frac{{\partial {p_{B - A}}}}{{\partial x}}} \right) + \frac{\partial }{{\partial y}}\left( {{f_2}\frac{{\partial {p_{B - A}}}}{{\partial y}}} \right) = f_3^\Delta ,(x,y) \in \Omega _{in}^B
\end{array} \right.}\]
where
\[{\left\{ \begin{array}{l}
f_3^\Delta  = 0,(x,y) \in {\Omega ^A} \cap (x,y) \in \Omega _{in}^A\\
f_3^\Delta  = {f_3},(x,y) \in \Omega _{in}^B \cap (x,y) \notin \Omega _{in}^A
\end{array} \right.}\]
According to the define of natural boundary condition:
\[{f_3} \ge 0,(x,y) \notin \Omega _{in}^A\]
we can get that:
\begin{equation}\label{p5_1}
f_3^\Delta  \ge 0
\end{equation}
via Hopf theorem:
\[\begin{array}{l}
{p_{B - A}}\left( {x,y} \right) \le 0\\
\exists (x,y) \in \Omega _{in}^B,{p_{B - A}}\left( {x,y} \right) < 0 \to \forall (x,y) \in \Omega _{in}^B,{p_{B - A}}\left( {x,y} \right) < 0
\end{array}\]
when ${p_{B - A}}\left( {x,y} \right) = 0$,the pressure obtained by $A,B$ boundary condition is actually the same.In other words, $A,B$  is actually a same solution.

Exclude this situation,the pressure is negative on every point out the boundary of A and in the boundary of B:
\[(x,y) \in {\Omega ^A} \cap (x,y) \in \Omega _{in}^B \to {p_{B - A}}\left( {x,y} \right) < 0 \to {p_B}\left( {x,y} \right) < 0 \to B \notin S\]

This conclusion is inconsistent with that the solution obtained through boundary $B$  is a positive solution ($B \in S$),so the original proposition is proved.

It should be pointed out that $\frac{{\partial p}}{{\partial {\Omega _n}}} = 0$ is indispensable for proof 5.1 as if $\frac{{\partial p}}{{\partial {\Omega _n}}} \ne 0$, $f_3$ is a negative impulse function on boundary ${\Omega _A}$ and formula \ref{p5_1} becomes invalid.

\begin{figure}[H] 
\centering 
\includegraphics[height=0.5\textwidth]{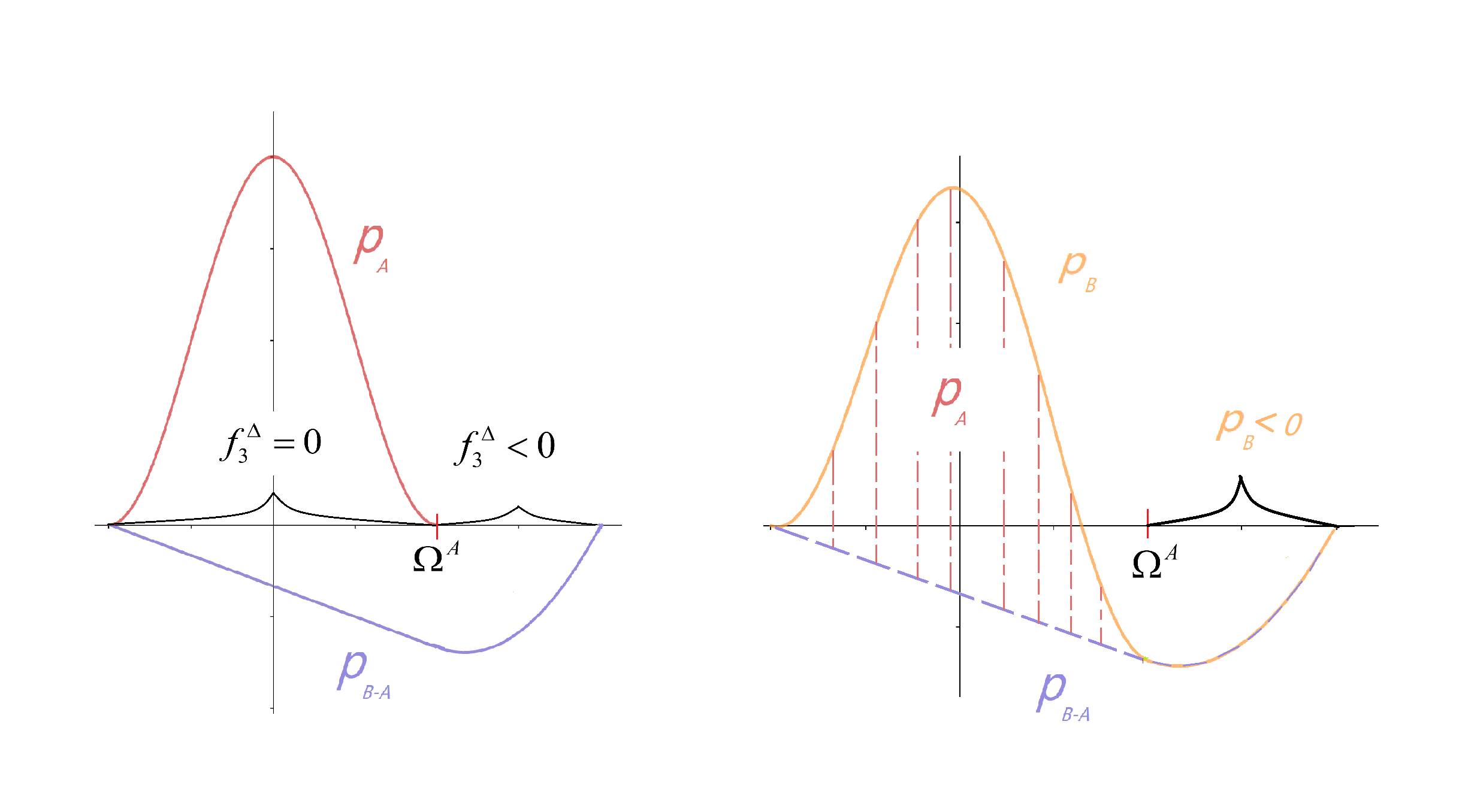} 
\caption{natural boundary condition must be maximum element in S} 
\label{Fig.NtoM} 
\end{figure}

Figure \ref{Fig.NtoM} is a simple diagram for this proof.

\noindent\textbf{proof 5.2.1}

Now let's prove the maximum element of $S$ must meet the natural boundary condition.

Suppose $A$ is a maximum element in $S$.According to the definition of $S$:
\[\left\{ \begin{array}{l}
{p_A}\left( {x,y} \right) > 0,(x,y) \in \Omega _{in}^A\\
{\left. {{p_A}} \right|_{{\Omega ^A}}} = 0
\end{array} \right.\]

If there is a point $P$ out boundary $\Omega_A$ where $f_3<0$, we can always find a neighborhood small enough where $f_3<0$ always correct.We can just add this neighborhood to $\Omega _{in}^A$ as boundary condition $A'$.Obviously, $A'$ is also in S and larger than $A$.This contradicts that A is the maximum element in S.\

Therefore we can get that:
\[\left\{ \begin{array}{l}
{p_A}\left( {x,y} \right) > 0,(x,y) \in \Omega _{in}^A\\
{\left. {{p_A}} \right|_{{\Omega ^A}}} = 0\\
{f_3} \ge 0,(x,y) \notin \Omega _{in}^A
\end{array} \right.\]

So there's only $\frac{{\partial p}}{{\partial {\Omega _n}}} = 0$ left to prove.As the pressure is always zero out the boundary $\Omega$,the $\frac{{\partial p}}{{\partial {\Omega _n}^ - }} = 0$ is always stays.Here, ${{\Omega _n}^ - } $ refers to the direction from outside to inside.

For the direction from inside to outside, as the pressure is larger than zero inside the boundary, $\frac{{\partial p}}{{\partial {\Omega _n}^ + }} > 0$ is impossible; we will proof $\frac{{\partial p}}{{\partial {\Omega _n}^ + }} < 0$ is also impossible in proof 5.2.2.

\noindent\textbf{proof 5.2.2}

Suppose $A$ is in set S,and there is a point ${\left( {{x_0},{y_0}} \right)}$ on $\Omega_A$ makes:
\[{\left. {\frac{{\partial {p_A}}}{{\partial {\Omega _n}^ + }}} \right|_{\left( {{x_0},{y_0}} \right)}} < 0\] 
We can take a neighborhood of $A$ on ${\left( {{x_0},{y_0}} \right)}$.And we can add this neighborhood to $A$ to create a new boundary condition $A'$. When the neighborhood is small enough,$A'$ must be an element of set S.As $A'>A$,$A$ can not be the maximum element.Therefore  $\frac{{\partial p}}{{\partial {\Omega _n}^ + }} < 0$ is impossible.

There is a simple proof idea. Small boundary changes lead to small changes in the first derivative of pressure.Therefore, the pressure derivative in the neighborhood remains negative, and the corresponding pressure in the neighborhood is positive.According to the linearity of the elliptic equation, the pressure solution is increased on the whole calculation area.So we can get that ${p'_A} > {p_A} > 0$.

\begin{figure}[H] 
\centering 
\includegraphics[width=0.7\textwidth]{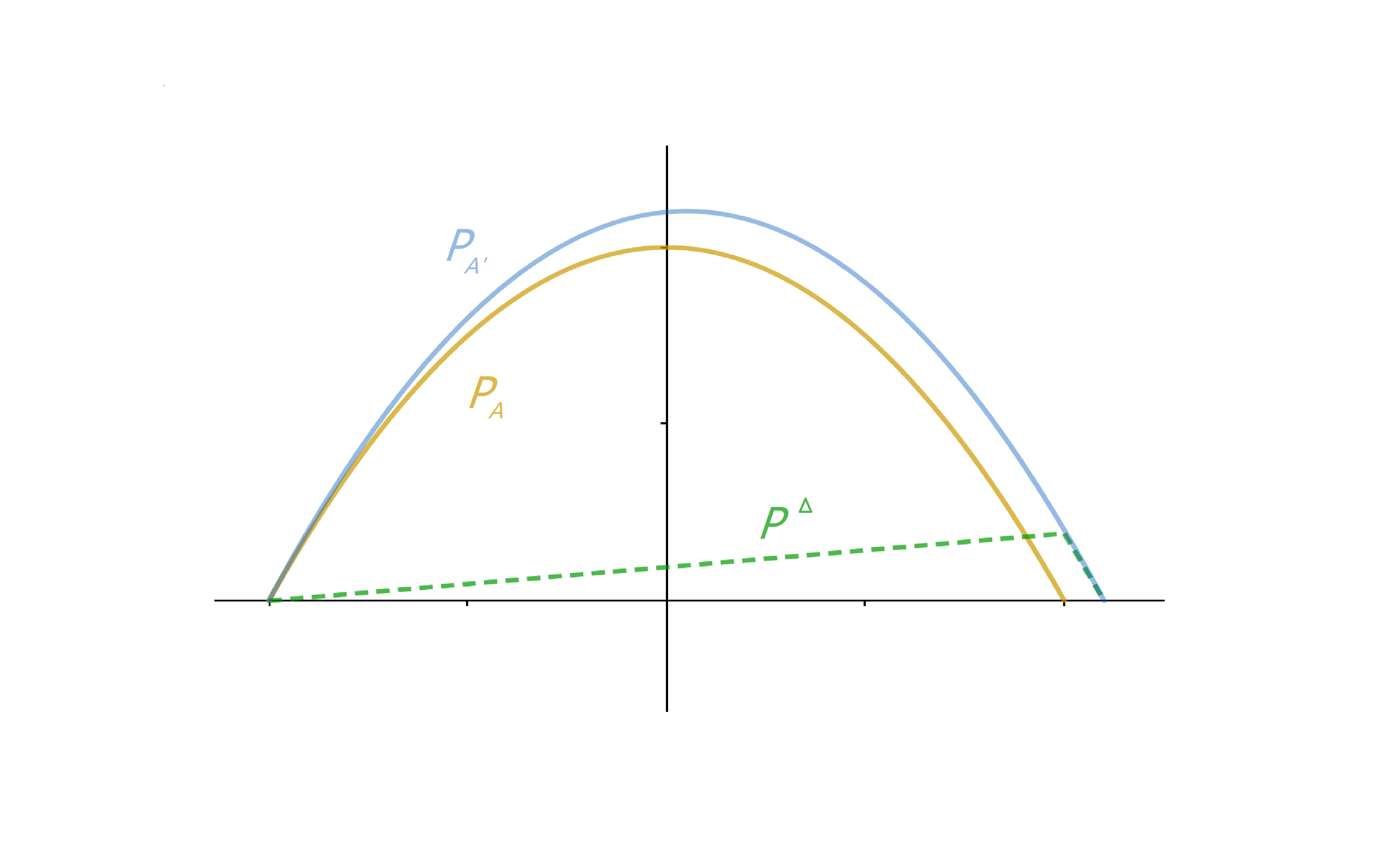} 
\caption{Neighborhood of A and A'} 
\label{Fig.NB} 
\end{figure}

Figure \ref{Fig.NB} is a sketch map. 

We can prove in more detail that there must be a neighborhood makes ${p'_A} > {p_A} > 0$.No mater rotating and translating the axes to make the point $\frac{{\partial p}}{{\partial {\Omega _n}^ + }} < 0$ on the point $(0,0)$,and $\Omega _n^ +  = {x^ + }$.  And scaling axes let in the area near the origin,to makes the x,y direction  have the same coefficientthe second order differential is 1.the elliptic equation has the following form:
\[\left\{ \begin{array}{l}
\frac{{{\partial ^2}p}}{{\partial {{x'}^2}}} + \frac{{{\partial ^2}p}}{{\partial {{y'}^2}}} = 1,\left( {x = 0,y = 0} \right)\\
\frac{{\partial p}}{{\partial x'}} =  - k,\left( {x = 0,y = 0} \right)
\end{array} \right.\]
where k is a positive real number.

Take an infinitesimal quantity tending to zero as $\lim t \to 0$ and we can give the neighborhood as :
\[\left\{ {\begin{array}{*{20}{l}}
{{{\left( {x - \frac{1}{3}k} \right)}^2} - \frac{1}{2}{y^2} - \frac{1}{9}{k^2} + \frac{1}{2}{t^2} > 0}\\
{0 < x < \frac{1}{3}k}
\end{array}} \right.\]
Add the above neighborhood to boundary $A$ to get boundary $A'$,assume that the solution of Reynolds equation under boundary $A$ is $p$,under boundary $A'$ is $p'$. We can get that: 
\[p' - p \ge {p_a} = \left\{ {\begin{array}{*{20}{l}}
{\max \left\{ {{{\left( {x - \frac{1}{3}k} \right)}^2} - \frac{1}{2}{y^2} - \frac{1}{9}{k^2} + \frac{1}{2}{t^2},{\rm{0}}} \right\},0 < x < \frac{1}{3}k}\\
{\max \left\{ {\frac{1}{2}{{\left( {x + \frac{1}{3}k} \right)}^2} - \frac{1}{2}{y^2} - \frac{1}{{18}}{k^2} + \frac{1}{2}{t^2},{\rm{0}}} \right\}, - \frac{1}{3}k < x < 0}
\end{array}} \right.\]
 It is not difficult to verify that the derivative of $p + {p_a}$ is zero in all area of $A'$ except on $\frac{1}{2}{\left( {x + \frac{1}{3}k} \right)^2} - \frac{1}{2}{y^2} - \frac{1}{{18}}{k^2} + \frac{1}{2}{t^2} = 0, - \frac{1}{3}k < x < 0$.Obviously, the $p$ here is 0 and $p \ge 0$ around this hyperbola. So the second derivative of ${p_a}$ is a negative impulse function.In other words the above formula stands.
 
 Therefore${\left. {\frac{{\partial {p_A}}}{{\partial {\Omega _n}^ + }}} \right|_{\left( {{x_0},{y_0}} \right)}} < 0$ contradiction with $A$ is the maximum element in $S$.

\[\neg \exists C \in S,C \ge A \to \left\{ \begin{array}{l}
\begin{array}{*{20}{l}}
{{p_A}\left( {x,y} \right) > 0,(x,y) \in \Omega _{in}^A}\\
{{{\left. {{p_A}} \right|}_{{\Omega ^A}}} = 0}\\
{{f_3} \ge 0,(x,y) \notin \Omega _{in}^A}
\end{array}\\
\left. {\frac{{\partial {p_A}}}{{\partial {\Omega _n}}}} \right| < 0
\end{array} \right.\]

We have proved the equivalence of the natural boundary condition and the maximum element in $S$.

\section{numerical proof}
\label{numerical_proof}
In practical application,the solution of Reynolds equation often use finite element analysis method to solve. We process the Reynolds equation with the discrete way to linear equations:

\begin{equation}\label{numerical}
\left\{ {\begin{array}{*{20}{l}}
\begin{array}{l}
k_{i,j}^1{p_{i + 1,j}} + k_{i,j}^2{p_{i - 1,j}} + k_{i,j}^3{p_{i,j + 1}} + k_{i,j}^4{p_{i,j - 1}} - {K_{i,j}}{p_{i,j}} = {D_{i,j}},\left( {i,j} \right) \in {\Omega _{in}}\\
{K_{i,j}} = \left( {k_{i,j}^1 + k_{i,j}^2 + k_{i,j}^3 + k_{i,j}^4} \right)
\end{array}\\
{{p_{i,j}} = 0,\left( {i,j} \right) \in \Omega }
\end{array}} \right.
\end{equation}

where ${k_{i,j}^1,k_{i,j}^2,k_{i,j}^3,k_{i,j}^4}$ and $D_{i,j}$ is known constant parameters, ${p_{i,j}}$ is pressure to be solved,$\Omega$ is given boundary of solution area.

Iterative method often be used to solve this question, and a well known method is used to make the solution meets the nature boundary condition as formula \ref{condition}: in each iteration of finite element analysis , we check all pressure of the finite element,if the pressure of an element is lower than 0, we set the pressure to 0:
\begin{enumerate}\label{itera}
\item get $p_{i,j}^{{T_{n + 0.5}}}$ by $p_{i,j}^{{T_{n}}}$
\item let $p_{i,j}^{{T_{n + 1}}} = \max \left( {p_{i,j}^{{T_{n + 0.5}}},0} \right)$
\item if $\forall \left( {p_{i,j}^{{T_{n + 1}}} - p_{i,j}^{{T_n}}} \right) < \delta$,solution end;else goto step 1
\end{enumerate}

In this way, the solution obtained after convergence will meet the natural boundary conditions.

We will prove that there is and only one set of pressure values $p_{i,j}$ to make the iteration error is 0:

\begin{equation}\label{z_error}
\forall \left( {p_{i,j}^{{T_{n + 1}}} - p_{i,j}^{{T_n}}} \right) = 0
\end{equation}

\noindent\textbf{Numerical form of hopf theorem}

Numerical form of hopf theorem: if in formula \ref{numerical}, only in one point $\left( {I,J} \right)$ ${D_{i,j}} \ne {\rm{0}}$ and ${D_{I,J}}{\rm{ < 0}}$,then ${p_{i,j}} > 0$:
\[\left\{ \begin{array}{l}
{D_{i,j}} < 0,\left( {i,j} \right) = \left( {I,J} \right)\\
{D_{i,j}} = 0,\left( {i,j} \right) \ne \left( {I,J} \right)
\end{array} \right. \to {p_{i,j}} > 0\]

This can be proved by following method: 

For ${\left( {i,j} \right) = \left( {I,J} \right)}$,
\[k_{i,j}^1{p_{i + 1,j}} + k_{i,j}^2{p_{i - 1,j}} + k_{i,j}^3{p_{i,j + 1}} + k_{i,j}^4{p_{i,j - 1}} - {K_{i,j}}{p_{i,j}} < 0\]
so ${p_{i,j}}$ cannot be the min value of ${p_{i + 1,j}},{p_{i - 1,j}},{p_{i,j + 1}},{p_{i,j - 1}}$.Naturally, ${p_{i,j}}$ cannot be the minimum value of all pressures.

For $\left( {i,j} \right) \ne \left( {I,J} \right),\left( {i,j} \right) \in \Omega$,
\[k_{i,j}^1{p_{i + 1,j}} + k_{i,j}^2{p_{i - 1,j}} + k_{i,j}^3{p_{i,j + 1}} + k_{i,j}^4{p_{i,j - 1}} - {K_{i,j}}{p_{i,j}} = 0\]
which means ${p_{i,j}}$ is the mean value of ${p_{i + 1,j}},{p_{i - 1,j}},{p_{i,j + 1}},{p_{i,j - 1}}$.
If ${p_{i,j}}$is the min value of all pressures, as all value of ${p_{i + 1,j}},{p_{i - 1,j}},{p_{i,j + 1}},{p_{i,j - 1}}$ can not be larger than ${p_{i,j}}$,so 
\[k_{i,j}^1{p_{i + 1,j}} + k_{i,j}^2{p_{i - 1,j}} + k_{i,j}^3{p_{i,j + 1}} + k_{i,j}^4{p_{i,j - 1}} - {K_{i,j}}{p_{i,j}} \ge 0\]
the equal sign is true when and only when ${p_{i + 1,j}} = {p_{i - 1,j}} = {p_{i,j + 1}} = {p_{i,j - 1}} = {p_{i,j}}$ is true.

And because of ${p_{i,j}}$is the min value of all pressures,we can get that ${p_{i + 1,j}}$ is also the min value.We can repeat the above proof to spread that the pressure is the same value on all $\Omega_{in}$.
This contradicts ${{D_{i,j}} < 0,\left( {i,j} \right) = \left( {I,J} \right)}$.Therefore ${p_{i,j}}$ also can not be the min value.
No matter whether $i,j$ is equal to $I,J$ or not, the pressure value of  ${p_{i,j}}$ cannot be the global minimum .So the min value must appear on the boundary which equals 0.

\noindent\textbf{The addition of positive solution}

It is not difficult to obtain the closeness of addition operation for the set of S using the method in Chapter 234.

The key point to find the solution is find in which points the pressure is zero,we can artificially divide the solution domain into two parts:

\begin{equation}\label{Solution_under_boundary}
\left\{ \begin{array}{l}
k_{i,j}^1{p_{i + 1,j}} + k_{i,j}^2{p_{i - 1,j}} + k_{i,j}^3{p_{i,j + 1}} + k_{i,j}^4{p_{i,j - 1}} - {K_{i,j}}{p_{i,j}} = {D_{i,j}},\left( {i,j} \right) \in \Omega _{in}^{(1)}\\
{p_{i,j}} = 0,\left( {i,j} \right) \notin \Omega _{in}^{(1)}
\end{array} \right.
\end{equation}
the part of ${p_{i,j}} = 0,\left( {i,j} \right) \notin \Omega _{in}^{(1)}$ can be used as the boundary condition to solve the pressure in $\left( {i,j} \right) \in \Omega _{in}^{(1)}$. If all pressure is not negative,we call this is a positive solution.

Suppose there are two boundary condition $\Omega _{in}^{(1)},\Omega _{in}^{(2)}$ are both positive solution,it is easy to prove that under a new boundary condition $\Omega _{in}^{(1 + 2)} = \Omega _{in}^{(1)} \cup \Omega _{in}^{(2)}$,the pressure is also not negative.Note $ p_{_{i,j}}^{\left( {1 + 2} \right)}$ 
is the solution under this boundary condition.Take $p_{_{i,j}}^{\max \left( {1,2} \right)} = \max \left( {p_{_{i,j}}^{(1)},p_{_{i,j}}^{(2)}} \right)$.

For any point, if 
\[p_{_{i,j}}^{\max \left( {1,2} \right)} = p_{_{i,j}}^{(1)}\]
as $p_{_{i,j}}^{\max \left( {1,2} \right)}$ is the larger of ${p_{_{i,j}}^{(1)},p_{_{i,j}}^{(2)}}$, we can know that $p_{_{i,j}}^{(1)}$ is greater than 0, so the following equation holds:
\[\begin{array}{l}
k_{i,j}^1p_{_{i + 1,j}}^{\max \left( {1,2} \right)} + k_{i,j}^2p_{_{i - 1,j}}^{\max \left( {1,2} \right)} + k_{i,j}^3p_{_{i,j + 1}}^{\max \left( {1,2} \right)} + k_{i,j}^4p_{_{i,j - 1}}^{\max \left( {1,2} \right)} - {K_{i,j}}p_{_{i,j}}^{\max \left( {1,2} \right)} - {D_{i,j}}\\
 \ge k_{i,j}^1p_{_{i + 1,j}}^{\left( 1 \right)} + k_{i,j}^2p_{_{i - 1,j}}^{\left( 1 \right)} + k_{i,j}^3p_{_{i,j + 1}}^{\left( 1 \right)} + k_{i,j}^4p_{_{i,j - 1}}^{\left( 1 \right)} - {K_{i,j}}p_{_{i,j}}^{\max \left( {1,2} \right)} - {D_{i,j}}\\
 = 0
\end{array}\]
According to linearity and Hopf's theorem:
\[p_{_{i,j}}^{\left( {1 + 2} \right)} \ge p_{_{i,j}}^{\max \left( {1,2} \right)} \ge 0\] 
It is not difficult to combine all the points appearing in all the positive solutions to form a new set, which is the set of all the points with pressure greater than 0 in the final solution:
\[\left\{ \begin{array}{l}
{\Omega ^{\left( {\rm{0}} \right)}}{\rm{ = }}{\Omega ^{\left( {\rm{1}} \right)}} \cup {\Omega ^{\left( {\rm{2}} \right)}} \cup {\Omega ^{\left( {\rm{3}} \right)}} \cup  \cdots  \cup {\Omega ^{\left( n \right)}}\\
p_{_{i,j}}^{(1)},p_{_{i,j}}^{(2)},p_{_{i,j}}^{(3)}, \cdots ,p_{_{i,j}}^{(n)} \ge 0
\end{array} \right.\]

The solution obtained by using boundary ${\Omega ^{\left( {\rm{0}} \right)}}$ as the boundary in formula \ref{Solution_under_boundary} is the final solution

\noindent\textbf{Iteration error is 0 only when the boundary is ${\Omega ^{\left( {\rm{0}} \right)}}$}

According to the iteration method \ref{itera} and the formula \ref{Solution_under_boundary}, we can get that:

\[\left\{ \begin{array}{l}
p_{_{i,j}}^{{T_{0.5}}} = p_{_{i,j}}^{\left( {\rm{0}} \right)} \ge 0,\left( {i,j} \right) \in \Omega _{in}^{({\rm{0}})}\\
p_{_{i,j}}^{{T_{0.5}}} = \frac{{k_{i,j}^1p_{_{i + 1,j}}^{\left( {\rm{0}} \right)} + k_{i,j}^2p_{_{i - 1,j}}^{\left( {\rm{0}} \right)} + k_{i,j}^3p_{_{i,j + 1}}^{\left( {\rm{0}} \right)} + k_{i,j}^4p_{_{i,j - 1}}^{\left( {\rm{0}} \right)} - {D_{i,j}}}}{{{K_{i,j}}}},\left( {i,j} \right) \notin \Omega _{in}^{({\rm{0}})}
\end{array} \right.\]

If there are $\left( {I,J} \right) \notin \Omega _{in}^{(0)}$, make $p_{_{i,j}}^{{T_{0.5}}} > 0$. We can get a new boundary as:
\[\Omega ' = {\Omega ^{\left( 0 \right)}} \cap \left\{ {\left( {i,j} \right) \notin \Omega _{in}^{(0)},p_{_{i,j}}^{{T_{0.5}}} > 0} \right\}\]
Suppose the new solution obtained by using boundary ${\Omega '}$ as the boundary in formula \ref{Solution_under_boundary} is $p',{p^\Delta } = p' - {p^{\left( {\rm{0}} \right)}}$,we can get that:
\[\begin{array}{l}
\left( {i,j} \right) \notin \Omega _{in}^{({\rm{0}})} \cap \left( {i,j} \right) \in {{\Omega '}_{in}}\\
 \to p_{_{i,j}}^{\left( 0 \right)} = 0,p_{_{i,j}}^{{T_{0.5}}} > 0\\
 \to p_{_{i,j}}^{\left( 0 \right)} < p_{_{i,j}}^{{T_{0.5}}} = \frac{{k_{i,j}^1p_{_{i + 1,j}}^{\left( 0 \right)} + k_{i,j}^2p_{_{i - 1,j}}^{\left( 0 \right)} + k_{i,j}^3p_{_{i,j + 1}}^{\left( 0 \right)} + k_{i,j}^4p_{_{i,j - 1}}^{\left( 0 \right)} - {D_{i,j}}}}{{{K_{i,j}}}}\\
 \to k_{i,j}^1p_{_{i + 1,j}}^{\left( 0 \right)} + k_{i,j}^2p_{_{i - 1,j}}^{\left( 0 \right)} + k_{i,j}^3p_{_{i,j + 1}}^{\left( 0 \right)} + k_{i,j}^4p_{_{i,j - 1}}^{\left( 0 \right)} - {K_{i,j}}p_{_{i,j}}^{\left( 0 \right)} > {D_{i,j}}
\end{array}\]
and 
\[k_{i,j}^1{{p'}_{i + 1,j}} + k_{i,j}^2{{p'}_{i - 1,j}} + k_{i,j}^3{{p'}_{i,j + 1}} + k_{i,j}^4{{p'}_{i,j - 1}} - {K_{i,j}}{{p'}_{i,j}} = {D_{i,j}},\left( {i,j} \right) \notin \Omega _{in}^{({\rm{0}})} \cap \left( {i,j} \right) \in {{\Omega '}_{in}}\]
we can get that:
\[\left\{ \begin{array}{l}
k_{i,j}^1p_{_{i + 1,j}}^\Delta  + k_{i,j}^2p_{_{i - 1,j}}^\Delta  + k_{i,j}^3p_{_{i,j + 1}}^\Delta  + k_{i,j}^4p_{_{i,j - 1}}^\Delta  - {K_{i,j}}p_{_{i,j}}^\Delta  = 0,\left( {i,j} \right) \in \Omega _{in}^{({\rm{0}})}\\
k_{i,j}^1p_{_{i + 1,j}}^\Delta  + k_{i,j}^2p_{_{i - 1,j}}^\Delta  + k_{i,j}^3p_{_{i,j + 1}}^\Delta  + k_{i,j}^4p_{_{i,j - 1}}^\Delta  - {K_{i,j}}p_{_{i,j}}^\Delta  < 0,\left( {i,j} \right) \notin \Omega _{in}^{({\rm{0}})} \cap \left( {i,j} \right) \in {{\Omega '}_{in}}\\
p_{_{i,j}}^\Delta  = 0,\left( {i,j} \right) \notin {{\Omega '}_{in}}
\end{array} \right.\]

according to the Hopf's theorem:
\[p_{_{i,j}}^\Delta  > 0\]

According to our above definition, ${\Omega '}$ should also be a subset of ${\Omega ^{\left( 0 \right)}}$.This contradicts $\Omega ' = {\Omega ^{\left( 0 \right)}} \cap \left\{ {\left( {i,j} \right) \notin \Omega _{in}^{(0)},p_{_{i,j}}^{{T_{0.5}}} > 0} \right\}$.So $p_{_{i,j}}^{{T_{0.5}}} \le 0$ and  $p_{_{i,j}}^{{T_{1}}} $ will be set to zero outside the boundary ${\Omega ^{\left( 0 \right)}}$.Which means:
\[p_{_{i,j}}^{{T_1}} = p_{_{i,j}}^{\left( 0 \right)}\]

That is to say, the iteration error of $p^{\left( 0 \right)}$ must be 0

Now we prove only ${\Omega ^{\left( 0 \right)}},p^{\left( 0 \right)}$ makes the iteration error is 0.
As ${\Omega ^{\left( 0 \right)}}$ contains all point in positive solutions,iteration error is 0 also means the solution is positive.Therefore any boundary makes iteration error is 0 is smaller than ${\Omega ^{\left( 0 \right)}}$. 
Suppose there is another boundary ${\Omega ^{\left( 1 \right)}} < {\Omega ^{\left( 0 \right)}}$ makes the iteration error is 0.

When $\left( {i,j} \right) \in \Omega _{in}^{(1)}$,we can get that:
\[\begin{array}{l}
k_{i,j}^1p_{i + 1,j}^{\left( 1 \right)} + k_{i,j}^2p_{i - 1,j}^{\left( 1 \right)} + k_{i,j}^3p_{i,j + 1}^{\left( 1 \right)} + k_{i,j}^4p_{i,j - 1}^{\left( 1 \right)} - {K_{i,j}}p_{i,j}^{\left( 1 \right)}\\
 = k_{i,j}^1p_{i + 1,j}^{\left( 0 \right)} + k_{i,j}^2p_{i - 1,j}^{\left( 0 \right)} + k_{i,j}^3p_{i,j + 1}^{\left( 0 \right)} + k_{i,j}^4p_{i,j - 1}^{\left( 0 \right)} - {K_{i,j}}p_{i,j}^{\left( 0 \right)}\\
 = {D_{i,j}}
\end{array}\]
When $\left( {i,j} \right) \notin \Omega _{in}^{(1)} \cap \left( {i,j} \right) \notin \Omega _{in}^{(0)}$,according to the iteration error is 0:
\[\begin{array}{*{20}{l}}
{p_{i,j}^{{T_1}} = p_{i,j}^{\left( 1 \right)} = 0}\\
{ \to p_{i,j}^{{T_{0.5}}} = \frac{{k_{i,j}^1p_{i + 1,j}^{\left( 1 \right)} + k_{i,j}^2p_{i - 1,j}^{\left( 1 \right)} + k_{i,j}^3p_{i,j + 1}^{\left( 1 \right)} + k_{i,j}^4p_{i,j - 1}^{\left( 1 \right)} - {D_{i,j}}}}{{{K_{i,j}}}} \le 0}\\
{ \to k_{i,j}^1p_{i + 1,j}^{\left( 1 \right)} + k_{i,j}^2p_{i - 1,j}^{\left( 1 \right)} + k_{i,j}^3p_{i,j + 1}^{\left( 1 \right)} + k_{i,j}^4p_{i,j - 1}^{\left( 1 \right)} - {K_{i,j}}p_{i,j}^{\left( 1 \right)} \le {D_{i,j}}}
\end{array}\]
therefore 
\[\begin{array}{l}
k_{i,j}^1p_{i + 1,j}^{\left( 0 \right)} + k_{i,j}^2p_{i - 1,j}^{\left( 0 \right)} + k_{i,j}^3p_{i,j + 1}^{\left( 0 \right)} + k_{i,j}^4p_{i,j - 1}^{\left( 0 \right)} - {K_{i,j}}p_{i,j}^{\left( 0 \right)}\\
 = {D_{i,j}} > \\
k_{i,j}^1p_{i + 1,j}^{\left( 1 \right)} + k_{i,j}^2p_{i - 1,j}^{\left( 1 \right)} + k_{i,j}^3p_{i,j + 1}^{\left( 1 \right)} + k_{i,j}^4p_{i,j - 1}^{\left( 1 \right)} - {K_{i,j}}p_{i,j}^{\left( 1 \right)}
\end{array}\]
according to the Hopf's theorem:
\[p_{i,j}^{\left( 0 \right)} < p_{i,j}^{\left( 1 \right)}\]
this obviously contradicts in the region outside ${\Omega ^{\left( 0 \right)}}$ and inside ${\Omega ^{\left( 1 \right)}}$, ${p ^{\left( 0 \right)}}$ is greater than 0 and ${p ^{\left( 1 \right)}}$ is equal to 0
\[\left( {i,j} \right) \notin \Omega _{in}^{(1)} \cap \left( {i,j} \right) \notin \Omega _{in}^{(0)} \to p_{i,j}^{\left( 0 \right)} > 0,p_{i,j}^{\left( 1 \right)} = 0\]

So we can get the result that Iteration error is 0 only when the boundary is ${\Omega ^{\left( {\rm{0}} \right)}}$ :
\[p_{i,j}^{_{{T_{n + 1}}}} = p_{i,j}^{{T_n}} \leftrightarrow \Omega  = {\Omega ^{\left( 0 \right)}}\]

\section{fast method for one dimension solution}
\label{fast}
For one dimension Reynolds equation,the discretized equation is a linear equation with tridiagonal matrix coefficients.And the solution get be got by algorithms with O (N) complexity.The only difficulty is the uncertainty of the boundary.According to our conclusions above, the boundary is always enlarged without shrinking, so we can gradually adjust the position of the boundary to get the final solution.

First, we discretize the Reynolds equation into:
\[{A_i}{p_{i - 1}} + {B_i}{p_{i + 1}} - \left( {{A_i} + {B_i}} \right){p_i} = {C_i}\]

According to the linearity of the equation we suppose that:
\[\left\{ {\begin{array}{*{20}{l}}
{{p_i} = {m_i}{p_1} + {n_i}{p_2} + {c_i}}\\
{{m_1} = 1,{n_1} = 0,{c_1} = 0}\\
{{m_2} = 0,{n_2} = 1,{c_2} = 0}
\end{array}} \right.\]

Therefore we can get that:
\[\begin{array}{l}
{A_i}\left( {{m_{i - 1}}{p_1} + {n_{i - 1}}{p_2} + {c_{i - 1}}} \right) + {B_i}\left( {{m_{i + 1}}{p_1} + {n_{i + 1}}{p_2} + {c_{i + 1}}} \right) - \left( {{A_i} + {B_i}} \right)\left( {{m_i}{p_1} + {n_i}{p_2} + {c_i}} \right) = {C_i}\\
 \to \left\{ \begin{array}{l}
{m_{i + 1}} = \frac{{\left( {{A_i} + {B_i}} \right){m_i} - {A_i}{m_{i - 1}}}}{{{B_i}}}\\
{n_{i + 1}} = \frac{{\left( {{A_i} + {B_i}} \right){n_i} - {A_i}{n_{i - 1}}}}{{{B_i}}}\\
{c_{i + 1}} = \frac{{\left( {{A_i} + {B_i}} \right){c_i} - {A_i}{c_{i - 1}} + {C_i}}}{{{B_i}}}
\end{array} \right.
\end{array}\]

Therefore, we can obtain ${m_i},{n_i},{c_i}$ before iterative solution.And according to the linearity of the equation, as long as we determine that the pressure value of any two points is 0 as the boundary condition, we can obtain the corresponding $p_1,p_2$. Then the pressure value at any point can be obtained.

Then we take all consecutive segments with ${C_i}<0$.Next,we use linked lists to store these segments, each linked list node stores the start and end points of the segment.

Next,we will iterate according to the following process:

\begin{enumerate}\label{itera}
\item Take the first node of the chain list as current node.
\item Assume that the pressure values at the starting point i and point j+1 of current node are 0.Where point j+1 is the next point of the ending point j.Calculate whether the pressure at j is less than 0. If it is greater than 0, go to step 3, otherwise go to step 4.
\item If there is a next node, and the starting point of the next node satisfy $i<j+1$ , merge the two nodes and take the new node as current node. Whether or not you merge, go back to step 2
\item If there is the next linked list node,take it as current node and go to step 2, otherwise go to step 5.
\item If this round of iteration moves the end point with a fixed starting point, the moving direction will be reversed to moves the starting point with a fixed end point and return to step 1, otherwise go to step 6.
\item If the starting and ending points do not move in the above two iterations (including the positive and negative directions), then go to step 7, otherwise go back to step 1.
\item According to the starting and ending points of each node, calculate the pressure value inside the node, and end the iteration.
\end{enumerate}

\begin{figure}[H] 
\centering 
\includegraphics[width=0.8\textwidth]{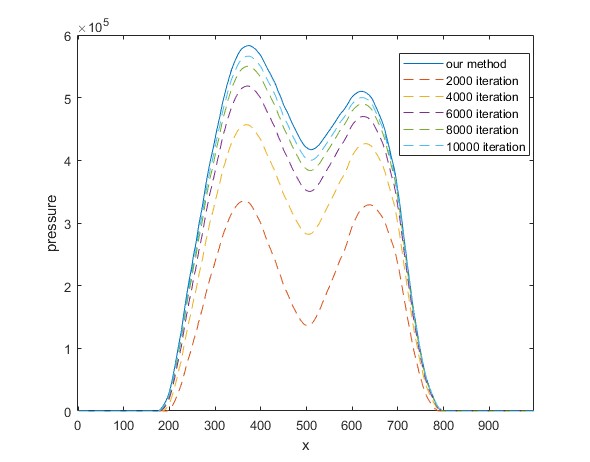} 
\caption{the test of our method and traditional iteration method} 
\label{Fig.valueTest} 
\end{figure}

We compare that our speed with traditional successive over relaxation (SOR) iteration method.Our method run 22.2 times faster than SOR. When we start iteration from only $5\%$ error, our running speed is 62.2 times that of SOR method. Our method complexity has reached O (N). When solving the Reynolds equation, the discretization of the Reynolds equation also requires O (N) computation. That is, the overall solution complexity will not be less than O (N). Therefore, our method has reached the theoretical upper limit of the complexity of solving the one-dimensional Reynolds equation.
The program is in cite \cite{github}.Figure \ref{Fig.valueTest} is a compare of our method result and traditional
iteration method.

\section{Conclusions}
\label{Conclusions}
We proved the existence and uniqueness of Reynolds equation under natural boundary conditions.This proof of the existence and uniqueness of the solution of the elliptic equation with natural boundary is of great significance to the lubrication problem and the phase change heat transfer problem.

Our work reveals the following two physical meanings for a wider range of fluid mechanics problems:

1.In the past, we believed that the physical basis of fluid dynamics equations was the motion of fluid microelements. Another explanation of fluid movement is that fluid naturally moves from high pressure area to low pressure area. For example, the propagation and expansion of shock wave generated by explosion. In our proof process, we can also see that the boundary line with pressure of 0 is essentially a greatest expansion. Therefore, it may be a new idea to regard the physical basis of fluid motion equation as the expansion of high pressure region. For example, the fluid dynamics equation is rewritten as the equation of the change of isobars with time and space, and its properties are studied.

2.In our research process, we found that in the one-dimensional case, we can transform the elliptic equation into a second-order differential equation through the expansion and contraction of coordinates. The function curve is obtained when the two endpoints are zero. Then, we find a line tangent to the curve from the bottom, and the region between the two tangent points is the boundary we are looking for. Therefore, we consider that this problem has an analytical solution in one-dimensional case. However, for high-dimensional cases, this coordinate transformation is not necessarily feasible. Inspired by this, we think that the topological properties of fluid equations after topological transformation may be a method to study the stability of fluid equations.

\bibliographystyle{model1-num-names}
\appendix

\bibliography{sample.bib}

\begin{thebibliography}{4}
\expandafter\ifx\csname natexlab\endcsname\relax\def\natexlab#1{#1}\fi
\providecommand{\url}[1]{\texttt{#1}}
\providecommand{\href}[2]{#2}
\providecommand{\path}[1]{#1}
\providecommand{\DOIprefix}{doi:}
\providecommand{\ArXivprefix}{arXiv:}
\providecommand{\URLprefix}{URL: }
\providecommand{\Pubmedprefix}{pmid:}
\providecommand{\doi}[1]{\href{http://dx.doi.org/#1}{\path{#1}}}
\providecommand{\Pubmed}[1]{\href{pmid:#1}{\path{#1}}}
\providecommand{\bibinfo}[2]{#2}
\ifx\xfnm\relax \def\xfnm[#1]{\unskip,\space#1}\fi
\bibitem[{Bayada et~al.(2007)Bayada, Talibi, and Hadi}]{bayada2007existence}
\bibinfo{author}{G.~Bayada}, \bibinfo{author}{M.~E.~A. Talibi},
  \bibinfo{author}{K.~A. Hadi},
\newblock \bibinfo{title}{Existence and uniqueness for a non-coercive
  lubrication problem},
\newblock \bibinfo{journal}{Journal of mathematical analysis and applications}
  \bibinfo{volume}{327} (\bibinfo{year}{2007}) \bibinfo{pages}{585--610}.
\bibitem[{Daquan and Tongchen(2010)}]{rapidmethod}
\bibinfo{author}{L.~Daquan}, \bibinfo{author}{M.~Tongchen},
\newblock \bibinfo{title}{The one-dimensional rapid algorithm for the
  generalized reynolds equations for journal bearings},
\newblock \bibinfo{journal}{Proceedings of the CSEE} \bibinfo{volume}{30}
  (\bibinfo{year}{2010}) \bibinfo{pages}{85--89}.
\bibitem[{Wen and Huang(2012)}]{wen2012principles}
\bibinfo{author}{S.~Wen}, \bibinfo{author}{P.~Huang},
  \bibinfo{title}{Principles of tribology}, \bibinfo{publisher}{John Wiley \&
  Sons}, \bibinfo{year}{2012}.
\bibitem[{git(1113)}]{github}
\bibinfo{title}{https://github.com/WangQunSJTU/solve-numerical-Reynolds-Equation},
  \bibinfo{year}{2022/11/13}.

\end{thebibliography}

\end{document}